\documentclass[12pt]{amsart}

\usepackage{amsmath}

\usepackage{amscd}

\usepackage{graphicx}

\usepackage{latexsym}

\usepackage{hyperref}

\usepackage{times}

\oddsidemargin .25in \theoremstyle{plain}

\def\a{\alpha}
\def\b{\beta}
\def\g{\gamma}
\def\s{\sigma}
\def\d{\delta}

\def\cp{\hbox{${\mathbb C} P^2$}}

\def\cpo{\hbox{${\mathbb C} P^1$}}

\def\cpb{\hbox{$\overline{{\mathbb C}P^2}$}}

\begin{document}

\title[]{On sections of elliptic fibrations}

\author{Mustafa Korkmaz}

\author{Burak Ozbagci}

\address{Department of Mathematics, Middle East Technical University,
 06531 Ankara, Turkey} \email{korkmaz@metu.edu.tr}

\address{School of Mathematics \\ Georgia Institute
of Technology \\  Atlanta  \\ Georgia and Department of Mathematics \\ Ko\c{c} University \\
Istanbul \\ Turkey}

\email{ bozbagci@math.gatech.edu and bozbagci@ku.edu.tr}


\date{\today}

\thanks{M.K. was partially supported by the Turkish Academy of Sciences.
B.O. was partially supported by the Turkish Academy of Sciences and by
the NSF Focused Research Grant FRG-024466.}

\begin{abstract}
We find a new relation among right-handed Dehn twists in the mapping
class group of a $k$-holed torus for  $4 \leq k \leq 9$. This
relation induces an elliptic Lefschetz pencil structure on the
four-manifold \cp $\#(9-k)$ \cpb $\;$ with $k$ base points and
twelve singular fibers. By blowing up the base points we get an
elliptic Lefschetz fibration on the complex elliptic surface
$E(1)=$\,\cp $\#9$ \cpb $\;\to S^2$ with twelve singular fibers and
$k$ disjoint sections. More importantly we can locate these $k$
sections in a Kirby diagram of the induced elliptic Lefschetz
fibration. The $n$-th power of our relation gives an explicit
description for $k$ disjoint sections of the induced elliptic
fibration on the complex elliptic surface $E(n) \to S^2$ for $n \geq
1$.

\end{abstract}

\maketitle

\section{Introduction}

It is well-known that two generic cubics $P$ and $Q$
in \cp$\;$ intersect each other in nine points $z_1,\ldots, z_9$.
By constructing
the corresponding \emph{pencil} of curves
$$\{ sP+tQ \;\vert\; [s:t] \in \cpo \} $$ one can define a map
$f : $ \cp $ - \{ z_1,\ldots, z_9\} \to $ \cpo. Blowing up \cp$\;$
at $\{z_1, \ldots , z_9\}$ one can extend $f$ to a Lefschetz
fibration $\pi : \; $\cp $\#9$ \cpb $\;\to \;$\cpo$\,$ whose generic
fiber is an elliptic curve. Our aim in this paper is to describe
this construction in terms of differential topology. It turns out
that many 4--manifold topologists were curious about such a
construction. (For instance this was posed explicitly as a question
in \cite{ga}).

Let $\Gamma_{g,k}^s$ denote the mapping class group of a compact connected orientable
genus $g$ surface with $k$ boundary components and $s$ marked points,
so that diffeomorphisms and isotopies of the surface are assumed to 
fix the marked points and the points on the boundary.
(We will drop $k$ if the surface is closed and drop $s$ if there are
no fixed points.) A product $\Pi_{i=1}^m t_i$ of right-handed Dehn
twists in $\Gamma_g$ provides a genus-$g$ Lefschetz fibration $X\to
D^2$ over the disk with closed fibers. If $\Pi_{i=1}^mt_i=1$ in
$\Gamma_g$ then the fibration closes up to a fibration over the
sphere $S^2$. A lift of the relation $\Pi_{i=1}^m t_i =1$ to
$\Gamma_g^k$ shows the existence of $k$ disjoint sections of the
induced Lefschetz fibration. The self-intersection of the $j$-th
section is $-n_j$ if $\Pi_{i=1}^m t_i =t_{\delta_1}^{n_1}\cdots
t_{\delta_k}^{n_k}$ in $\Gamma_{g,k}\;$, for some positive integers
$n_1, \ldots, n_k$,
 where $t_{\delta_i}$'s are right-handed Dehn twists along circles
parallel to the
  boundary components of the surface at hand (cf. \cite{ekkos}).

On the other hand, an expression $\Pi_{i=1}^m t_i=t_{\delta_1}
\cdots t_{\delta_k}$ in
  $\Gamma_{g,k}\;$ naturally describes a Lefschetz pencil: The relation determines a
Lefschetz fibration with $k$ disjoint sections, where each section
has self-intersection $-1$, and after blowing these sections down we
get a Lefschetz pencil (cf. \cite{ga}). Conversely, by blowing up
the base locus of a Lefschetz pencil we arrive to a Lefschetz
fibration which can be captured (together with the exceptional
divisors of the blow-ups, which are all sections now) by a relation
of the above type.

In this paper we find relations of the form $\Pi_{i=1}^{12}
t_i=t_{\delta_1}\cdots t_{\delta_k}$ in $\Gamma_{1,k}$ for $4 \leq k
\leq 9$, generalizing the well-known cases $k=1,2,3$. A relation of
this type naturally induces a Lefschetz pencil and by blowing up we
get an elliptic Lefschetz fibration with $k$ sections. Moreover by
taking the $n$-th power of our relation (for $n\geq 1$) we get
$$(\Pi_{i=1}^{12} t_i)^n =t_{\delta_1}^n \cdots t_{\delta_k}^n \in \Gamma_{1,k}$$
for $4 \leq k \leq 9$. Once again this relation induces an elliptic
Lefschetz fibration $E(n) \to S^2$ with $12n$ singular fibers and
$k$ disjoint sections, where the self-intersection of each section
is equal to $-n$.

The reader is advised to turn to \cite{d},  \cite{gs} and
\cite{ozst} for background material on Lefschetz fibrations and
pencils. To simplify the notation in the rest of the paper we will
denote a right-handed Dehn twist along a curve along $\alpha$ also
by $\alpha$. A left-handed Dehn twist along $\alpha$ will be denoted
by $\overline{\alpha}$. We will multiply the Dehn twists from right
to left, i.e., $\b\a$ means that we first apply $\a$ then $\b$.

\vspace{1ex}

\noindent {\bf {Acknowledgement}}: The authors would like to thank
John Etnyre, David Gay and Andras Stipsicz for their interest in
this work.

\section{Lantern relation for the four-holed sphere}\label{lant}
\begin{figure}[ht]

  \begin{center}

     \includegraphics{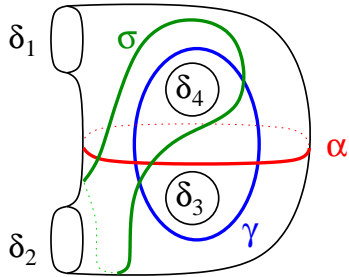}

   \caption{Four-holed sphere with boundary $\{ \d_1,\d_2,\d_3,\d_4 \}. $}

  \label{lantern}

    \end{center}

  \end{figure}

Consider the four-holed sphere depicted in Figure~\ref{lantern}.
Then there is the relation
\begin{eqnarray*}
\d_1\d_2\d_3 \d_4 = \g \s \a
\end{eqnarray*}
in $\Gamma_{0,4}$ which was discovered by Dehn~\cite{de}. It was
rediscovered by Johnson~\cite{jo} and named as the lantern relation.
We will freely use this relation in any subsurface (of another
surface) which is homeomorphic to a sphere with four holes.  The
particular depiction of the lantern relation on the four-holed
sphere in Figure~\ref{lantern} will be convenient in the subsequent
discussions.

\section{Relations on torus with holes}
In this section we will generalize the well known one-holed torus
relation to a relation on the $k$-holed torus for $2 \leq k \leq 9$.
We will give all the details in each case since the relation for
$(k+1)$-holes is derived by using the relation for $k$-holes for $1
\leq k \leq 8$. The relations in the cases $k=2,3$ are also known,
but we compute them anyway for the sake of completeness and to show
our method.

We note that if two circles are disjoint, then the corresponding Dehn twists
commute. Also, if two circles $\alpha$ and $\beta$ intersect transversely at one
point, then the corresponding Dehn twists satisfy the braid relation;
$\alpha\beta \alpha = \beta\alpha\beta$. We use these facts in our computations below.

\subsection{One-holed torus}\label{one-holedtorus}
If $\a$ and $\b$ are two circles on a torus with one boundary
$\d_1$ which intersect each other transversely at one point, then
there is the relation
\begin{eqnarray*}
(\a \b)^6=\d_1.
\end{eqnarray*}
We call it the one-holed torus relation. It turns out that this relation was known to
Dehn~\cite{de} in a slightly different form.

\subsection{Two-holed torus}
\begin{figure}[ht]

  \begin{center}

     \includegraphics{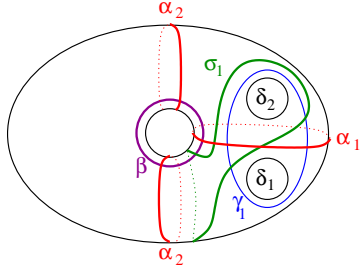}

     \caption{Two-holed torus with boundary $\{ \d_1,\d_2\}. $}

  \label{2-star}

    \end{center}

  \end{figure}

Consider the two-holed torus depicted in Figure~\ref{2-star}. By
the lantern relation, we have
$$\a_2^2\d_1\d_2=\g_1\s_1\a_1 .$$
One-holed torus relation is
$$\g_1=(\a_2 \b)^6.$$
Note that $\s_1=\overline{\b}\overline{\a}_2\overline{\a}_2 \a_1 \b
 \overline{\a}_1\a_2\a_2\b$.
Then
\begin{eqnarray*}
\d_1\d_2&=& \overline{\a}_2\,\overline{\a}_2 \g_1\s_1  \a_1\\
&=& \overline{\a}_2\,\overline{\a}_2 (\a_2  \b\a_2 \b \a_2 \b \a_2 \b \a_2 \b
\a_2 \b ) (\overline{\b}\overline{\a}_2\overline{\a}_2 \a_1 \b
 \overline{\a}_1\a_2\a_2\b) \a_1\\
&=& \overline{\a}_2\,\overline{\a}_2 \a_2 \a_2 \b \a_2 \a_2 \b \a_2 \b \a_2 \b \overline{\a}_2 \a_1 \b
 \overline{\a}_1\a_2\a_2\b \a_1\\
&=& \b \a_2 \a_2 \b \a_2 \a_2 \b   \a_1 \b \overline{\a}_1\a_2\a_2\b \a_1\\
&=& \b \a_2 \a_2 \b \a_2 \a_2  \a_1 \b \a_2\a_2\b \a_1.
\end{eqnarray*}
It follows that $\d_1\d_2=(\a_1 \a_2 \b)^4$.

\subsection{Three-holed torus} \label{three-star}
\begin{figure}[ht]

  \begin{center}

     \includegraphics{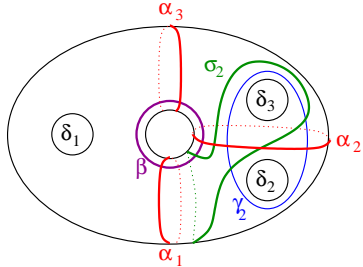}

   \caption{Three-holed torus with boundary $\{ \d_1,\d_2,\d_3\}. $}

  \label{3-star}

    \end{center}

  \end{figure}

Consider the three-holed torus depicted in Figure~\ref{3-star}. By
the lantern relation,
$$\a_3 \a_1\d_2\d_3=\g_2\s_2\a_2$$
and by the two-holed torus relation,
$$\d_1\g_2=(\a_1\a_3\b)^4.$$
Note that $\s_2=\overline{\b}\overline{\a}_1\overline{\a}_3 \a_2
\b \overline{\a}_2\a_3\a_1\b$.
Then
\begin{eqnarray*}
\d_1\d_2\d_3 &=& \overline{\a}_1\,\overline{\a}_3 \d_1 \g_2\s_2  \a_2\\
&=& \overline{\a}_1\,\overline{\a}_3
(\a_1 \a_3 \b\a_1 \a_3 \b\a_1 \a_3 \b\a_1 \a_3 \b ) (\overline{\b}\overline{\a}_1 \overline{\a}_3 \a_2 \b
 \overline{\a}_2\a_3\a_1\b) \a_2\\
&=&  \b\a_1 \a_3 \b\a_1 \a_3 \b  \a_2 \b \overline{\a}_2\a_3\a_1\b  \a_2\\
&=&  \b\a_1 \a_3 \b\a_1 \a_3  \a_2 \b  \a_3\a_1\b  \a_2.
\end{eqnarray*}
It follows that $\d_1\d_2\d_3 =  (\a_1\a_2\a_3 \b)^3$. This relation
was called the star relation in~\cite{ge}.

\subsection{Four-holed torus}\label{four-star}

\begin{figure}[ht]

  \begin{center}

     \includegraphics{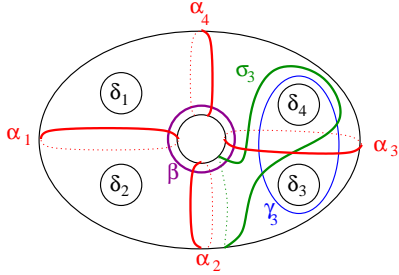}

   \caption{Four-holed torus with boundary $\{ \d_1,\d_2,\d_3, \d_4\}. $}

  \label{4-star}

    \end{center}

  \end{figure}

The lantern relation for the sphere with boundary
$\{\a_4,\a_2,\d_3,\d_4\}$ in Figure~\ref{4-star} is $$\a_4\a_2
\d_3\d_4= \g_3\s_3\a_3.$$
The relation on the three-holed
torus with boundary $\{\d_1,\d_2,\g_3 \}$ given in Section~\ref{three-star} is
$$\d_1\d_2\g_3 = (\a_1\a_2\a_4\b)^3.$$ Here we identify the curves  $(\a_1,\a_2,\a_3)$
in Figure~\ref{3-star} by the curves $(\a_1,\a_2,\a_4)$ in Figure~\ref{4-star}.  Combining we get
\begin{eqnarray*}
\d_1\d_2\d_3\d_4 &=& \d_1\d_2 \,
\overline{\a}_2\,\overline{\a}_4\g_3\s_3 \a_3\\ &=&
\overline{\a}_2\,\overline{\a}_4\d_1\d_2 \g_3 \s_3 \a_3\\ &=&
\overline{\a}_2\,\overline{\a}_4(\a_1\a_2\a_4\b)^3\s_3 \a_3\\ &=&
\a_1\b(\a_1\a_2\a_4\b)^2\s_3\a_3
\\ &=& (\a_1\a_2\a_4\b)^2\s_3\a_3\a_1\b.
\end{eqnarray*}

\noindent{\bf Remark.} Although we will not need it in the rest of the paper, by
plugging in
$$\s_3=\overline{\b}\overline{\a}_4\overline{\a}_2 \a_3\b
\overline{\a}_3\a_2\a_4\b,$$ it is easy to see that this relation
may also be written in a more symmetric form as
\begin{eqnarray*}
\d_1\d_2\d_3\d_4 = \left(\a_1\a_3 \b\a_2\a_4\b)^2. \right.
\end{eqnarray*}

\subsection{Five-holed torus}\label{five-star}

\begin{figure}[ht]

  \begin{center}

     \includegraphics{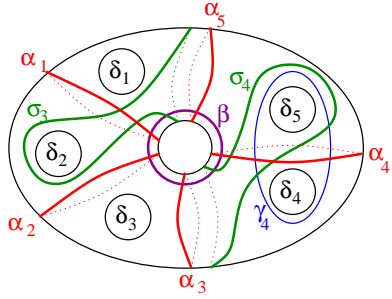}

   \caption{Five-holed torus with boundary $\{ \d_1,\d_2,\d_3,\d_4, \d_5\}. $}

  \label{5-star}

    \end{center}

  \end{figure}

The lantern relation for the sphere with boundary
$\{\a_5,\a_3,\d_4,\d_5\}$ in Figure~\ref{5-star} is
$$\a_5\a_3\d_4\d_5= \g_4\s_4\a_4.$$
The relation on the four-holed torus with boundary
$\{\d_1,\d_2,\d_3,\g_4\}$ given in Section~\ref{four-star} is
$$\d_1\d_2\d_3\g_4 = (\a_3\a_5\a_2\b)^2\s_3\a_1\a_3\b.$$ Here we
identify the curves $(\a_1,\a_2,\a_3,\a_4)$ in Figure~\ref{4-star}
with the curves $(\a_3,\a_5,\a_1,\a_2)$ in Figure~\ref{5-star}.
Combining we get

\begin{eqnarray*}
\d_1\d_2\d_3\d_4\d_5 &=& \overline{\a}_3\,\overline{\a}_5
\d_1\d_2\d_3 \g_4 \s_4\a_4 \\ &=&
\overline{\a}_3\,\overline{\a}_5(\a_3\a_5\a_2\b)^2\s_3\a_1\a_3\b\s_4\a_4
\\ &=& \a_2\b\a_3\a_5\a_2\b\s_3\a_1\a_3\b\s_4 \a_4\\ &=&
\a_2\a_3\a_5\b\s_3\a_1\a_3\b\s_4\a_4\a_2\b.
\end{eqnarray*}

\subsection{Six-holed torus} \label{six-star} 
\begin{figure}[ht]

  \begin{center}

     \includegraphics{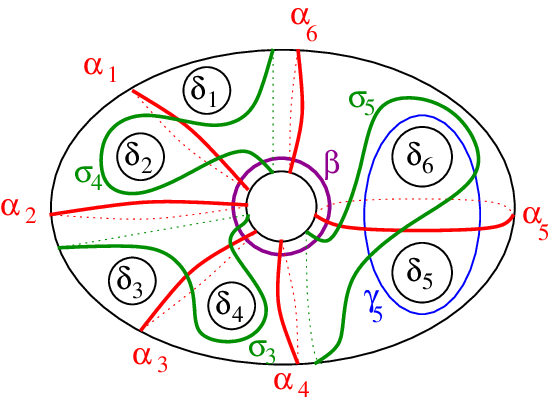}

   \caption{Six-holed torus with boundary $\{ \d_1,\d_2,\ldots, \d_6\}. $}

  \label{6-star}

    \end{center}

  \end{figure}

The lantern relation for the
sphere with boundary $\{\a_6, \a_4, \d_5,\d_6\}$ in
Figure~\ref{6-star} is
$$\a_6\a_4\d_5\d_6 = \g_5\s_5\a_5.$$
The relation for the five-holed torus with boundary
$\{\d_1,\d_2,\d_3,\d_4,\g_5\}$ given in Section~\ref{five-star} is
$$\d_1\d_2\d_3\d_4\g_5 = \a_4\a_6\a_2\b\s_3\a_3\a_6\b\s_4\a_1\a_4\b.$$
We identify the curves $(\a_1,\a_2,\a_3,\a_4,\a_5) $ in
Figure~\ref{5-star} with the curves $(\a_3,\a_4,\a_6,\a_1,\a_2)$ in
Figure~\ref{6-star}. Combining we get

\begin{eqnarray*}
\d_1\d_2\d_3\d_4\d_5\d_6 &=& \overline{\a}_4\,\overline{\a}_6
\d_1\d_2\d_3\d_4  \g_5 \s_5\a_5
\\ &=&\overline{\a}_4\,\overline{\a}_6\a_4\a_6\a_2\b\s_3\a_3\a_6\b\s_4\a_1\a_4\b\s_5\a_5
\\ &=& \a_2\b\s_3\a_3\a_6\b\s_4\a_1\a_4\b\s_5 \a_5\\ &=&
\b_2\a_2\s_3\a_3\a_6\b\s_4\a_1\a_4\b\s_5\a_5 \\&=&
\a_2\a_3\a_6\b\s_4\a_1\a_4\b\s_5\a_5\b_2\s_3,
\end{eqnarray*}
where $\b_2= \a_2 \b \overline{\a}_2$.

\subsection{Seven-holed torus} \label{seven-star} 
\begin{figure}[ht]

  \begin{center}

     \includegraphics{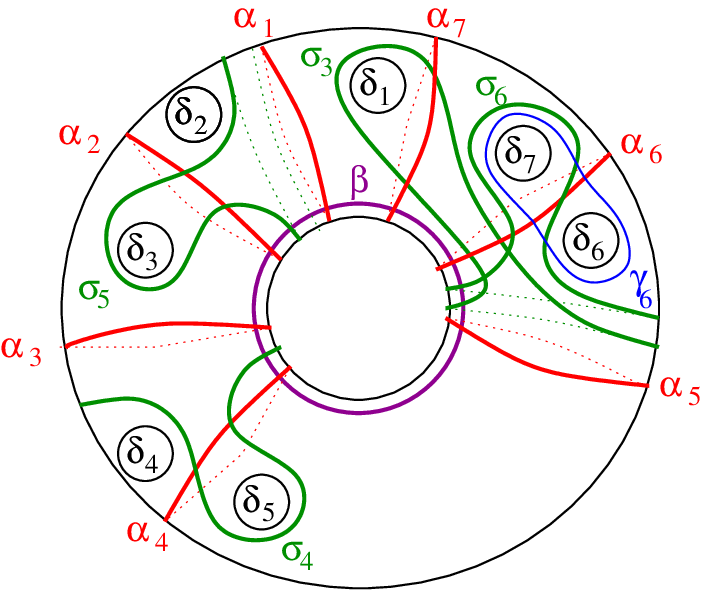}

   \caption{Seven-holed torus with boundary $\{ \d_1,\d_2,\ldots,\d_7\}. $}

  \label{7-star}

    \end{center}

  \end{figure}

The lantern relation for the
sphere with boundary $\{\a_7,\a_5,\d_6,\d_7\}$ in
Figure~\ref{7-star} is
$$\a_7\a_5\d_6\d_7= \g_6\s_6\a_6.$$
The relation on the six-holed torus with boundary
$\{\d_1,\d_2,\d_3,\d_4,\d_5,\g_6\}$ given in Section~\ref{six-star} is
$$\d_1\d_2\d_3\d_4\d_5\g_6 =
\a_5\a_7\a_3\b\s_4\a_4\a_1\b\s_5\a_2\b_5\s_3,$$ where we use the
identification $(\a_1,\a_2,\a_3,\a_4,\a_5,\a_6) \to
(\a_4,\a_5,\a_7,\a_1,\a_2,\a_3)$ to go from Figure~\ref{6-star} to
Figure~\ref{7-star}. Combining we get

\begin{eqnarray*}
\d_1\d_2\d_3\d_4\d_5\d_6\d_7 &=& \overline{\a}_5\,
\overline{\a}_7\d_1\d_2\d_3\d_4\d_5 \g_6\s_6\a_6
\\ &=& \overline{\a}_5 \, \overline{\a}_7\a_5\a_7\a_3\b\s_4\a_4\a_1\b\s_5\a_2\b_5\s_3\s_6\a_6
\\ &=& \a_3\b\s_4\a_4\a_1\b\s_5\a_2\b_5\s_3\s_6\a_6\\ &=&
\b_3\a_3\s_4\a_4\a_1\b\s_5\a_2\b_5\s_3\s_6\a_6\\&=&
\b_3\s_4\a_3\a_4\a_1\b\s_5\a_2\b_5\s_3\s_6\a_6\\&=&
\a_3\a_4\a_1\b\s_5\a_2\b_5\s_3\s_6\a_6\b_3\s_4,
\end{eqnarray*}
where $\b_3= \a_3 \b \overline{\a}_3$ and $\b_5= \a_5 \b \overline{\a}_5$.

\subsection{Eight-holed torus} \label{eight-star} 
\begin{figure}[ht]

  \begin{center}

     \includegraphics{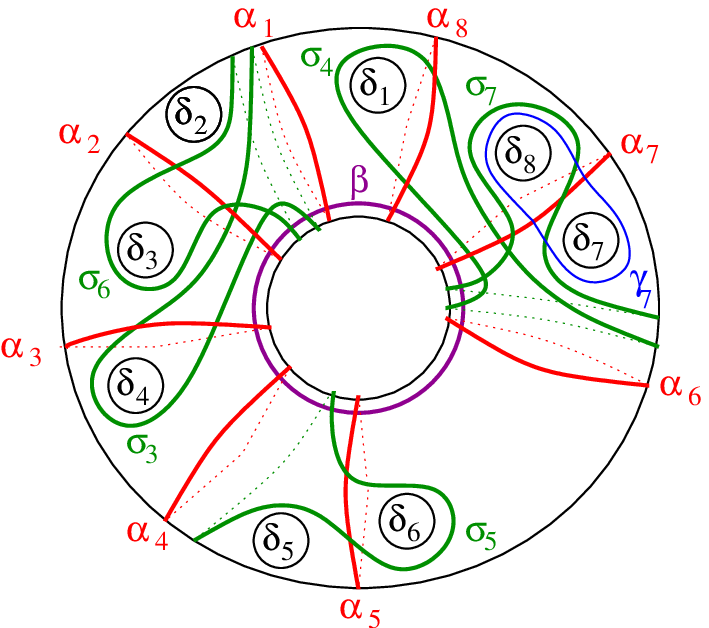}

   \caption{Eight-holed torus with boundary $\{ \d_1,\d_2, \ldots,\d_8\}. $}

  \label{8-star}

    \end{center}

  \end{figure}

The lantern relation for the
sphere with boundary $\{\a_8,\a_6,\d_7,\d_8\}$ in
Figure~\ref{8-star} is
$$\a_8\a_6\d_7\d_8 = \g_7\s_7\a_7.$$
The relation on the seven-holed torus with boundary
$\{\d_1,\d_2,\d_3,\d_4,\d_5,\d_6,\g_7\}$  given in Section~\ref{seven-star} is
$$\d_1\d_2\d_3\d_4\d_5\d_6\g_7=\a_6\a_8\a_4\b\s_5\a_5\b_1\s_3\s_6\a_2\b_6\s_4,$$
where we use the identification
$(\a_1,\a_2,\a_3,\a_4,\a_5,\a_6,\a_7) \to
(\a_4,\a_5,\a_6,\a_8,\a_1,\a_2,\a_3)$ to go from Figure~\ref{7-star}
to Figure~\ref{8-star}. Combining we get

\begin{eqnarray*}
\d_1\d_2\d_3\d_4\d_5\d_6\d_7\d_8 &=&
\overline{\a}_6\, \overline{\a}_8\d_1\d_2\d_3\d_4\d_5\d_6 \g_7\s_7\a_7\\
&=& \overline{\a}_6\, \overline{\a}_8\a_6\a_8\a_4\b\s_5\a_5\b_1\s_3\s_7\a_2\b_6\s_4\s_7\a_7\\
&=& \a_4\b\s_5\a_5\b_1\s_3\s_6\a_2\b_6\s_4\s_7\a_7\\ &=&
\b_4\a_4\s_5\a_5\b_1\s_3\s_6\a_2\b_6\s_4\s_7\a_7\\ &=&
\b_4\s_5\a_4\a_5\b_1\s_3\s_6\a_2\b_6\s_4\s_7\a_7\\ &=&
\a_4\a_5\b_1\s_3\s_6\a_2\b_6\s_4\s_7\a_7\b_4\s_5,
\end{eqnarray*}
where $\b_1= \a_1 \b \overline{\a}_1$, $\b_4= \a_4 \b \overline{\a}_4$ and
$\b_6= \a_6 \b \overline{\a}_6$.

\subsection{Nine-holed torus} \label{nine-star} 
\begin{figure}[ht]

  \begin{center}

     \includegraphics{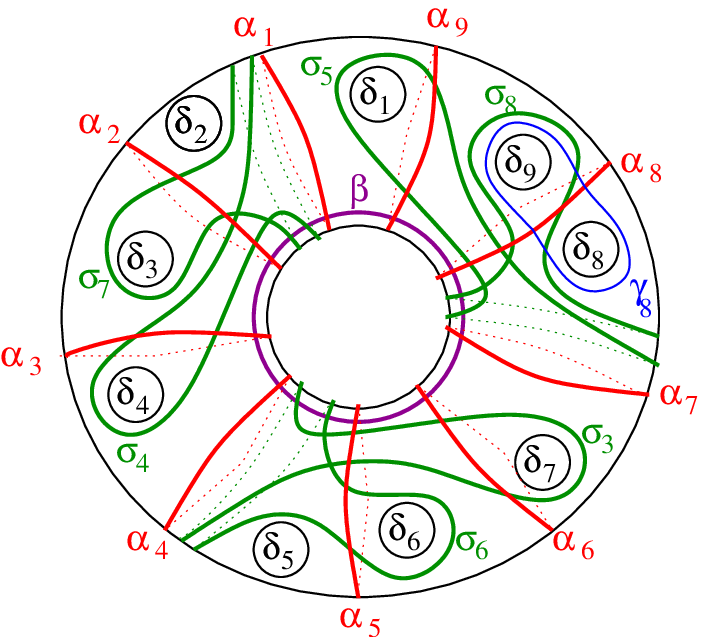}

   \caption{Nine-holed torus with boundary $\{ \d_1,\d_2,\ldots, \d_9\}. $}

  \label{9-star}

    \end{center}

  \end{figure}

The lantern relation for the
sphere with boundary $\{\a_9, \a_7,\d_8,\d_9\}$ in
Figure~\ref{9-star} is
$$\a_9\a_7\d_8\d_9= \g_8\s_8\a_8.$$
The relation on the eight-holed torus with boundary
$\{\d_1,\d_2,\d_3,\d_4,\d_5,\d_6,\d_7 \g_8\}$  given in Section~\ref{eight-star} is
$$\d_1\d_2\d_3\d_4\d_5\d_6\d_7\g_8=\a_7\a_9\b_4\s_3\s_6\a_5\b_1
\s_4\s_7\a_2\b_7\s_5,$$ where we identify
$(\a_1,\a_2,\a_3,\a_4,\a_5,\a_6,\a_7,\a_8)$ with
$(\a_4,\a_5,\a_6,\a_7,\a_9,\a_1,\a_2,\a_3)$ to go from
Figure~\ref{8-star} to Figure~\ref{9-star}. Combining we get

\begin{eqnarray*}
\d_1\d_2\d_3\d_4\d_5\d_6\d_7\d_8\d_9 &=&
\overline{\a}_7\,\overline{\a}_9\d_1\d_2\d_3\d_4\d_5\d_6\d_7\g_8\s_8\a_8\\
&=& \overline{\a}_7\,\overline{\a}_9\a_7\a_9\b_4\s_3\s_6\a_5\b_1
\s_4\s_7\a_2\b_7\s_5\s_8\a_8\\ &=& \b_4\s_3\s_6\a_5\b_1
\s_4\s_7\a_2\b_7\s_5\s_8\a_8,
\end{eqnarray*}
where $\b_1= \a_1 \b \overline{\a}_1$, $\b_4= \a_4 \b \overline{\a}_4$ and
$\b_7= \a_7 \b \overline{\a}_7$.

\bigskip
\noindent{\bf Remark.} The curious reader might wonder why we stopped
at $k=9$. First of all our process will not allow us to go any
further. In fact, there is a good reason for that: An elliptic
fibration $E(1) \to S^2$ admits at most nine disjoint sections (all
with negative self-intersections). So we conclude that there is no
such relation for $k$-holed torus with $k \geq 10$.

\section{Sections of the elliptic fibrations}
\label{sections}

\begin{figure}[ht]

  \begin{center}

     \includegraphics[width=6cm]{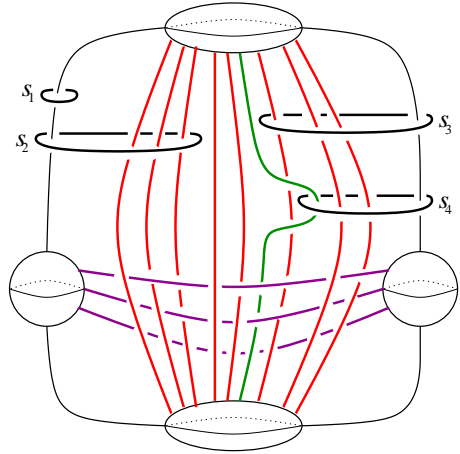}

   \caption{An elliptic Lefschetz fibration $E(1) \to S^2$ with four disjoint sections. }

  \label{4-section}

    \end{center}

  \end{figure}

First we consider the case $k=4$. The relation
$$\d_1\d_2\d_3\d_4 = (\a_1\a_2\a_4\b)^2\s_3\a_3\a_1\b $$
 in $\Gamma_{1,4}$ we derived in Section~\ref{four-star} induces the word
$(\a^3\b)^3=1$ in $\Gamma_1$ which gives us an elliptic Lefschetz
fibration on the elliptic surface $E(1)=$\,\cp $\#9$ \cpb. To draw a
Kirby diagram (\cite{gs}) of this elliptic fibration we start with a
0-handle, attach two 1-handles (see Figure~\ref{4-section}) and
attach a 2-handle which yields $D^2 \times T^2$.  A torus fiber of
the trivial fibration $D^2 \times T^2 \to D^2$ can be viewed in
Figure~\ref{4-section} as follows: Take the obvious disk on the
page, attach two 2-dimensional 1-handles (going through two
4-dimensional 1-handles) and cap off by a 2-dimensional disk. Then
we draw the curves which appear in the monodromy of the elliptic
fibration on parallel copies of this fiber. Notice that these curves
are the attaching curves of some 2-handles. Once we attach all
twelve of these 2-handles with framing one less than the page
framing we get an elliptic Lefschetz fibration over $D^2$ with
twelve singular fibers, which then can be closed off to an elliptic
Lefschetz fibration over $S^2$. We depicted the four disjoint
sections $s_1, s_2,s_3, s_4$ of the induced fibration in
Figure~\ref{4-section}. (Imagine replacing $s_i$'s in
Figure~\ref{4-section} by holes where they intersect the page and
embed the curves in Figure~\ref{4-star} into distinct fibers.) For
each $i=1,2,3,4$, the curve $s_i$ bounds two disks---one in the
neighborhood of a regular fiber, one outside of that
neighborhood---which in turn gives a section of the elliptic
Lefschetz fibration when we glue these two disks along their common
boundary $s_i$.

Similarly we can draw the Kirby diagrams corresponding to the
relations we derived for $k=5,6,\ldots,9$, and explicitly indicate
the locations of the $k$ disjoint sections of $E(1) \to S^2$ in
these diagrams. We skip the cases $k=5, \ldots, 8$ and jump to the
case $k=9$. The relation
$$\d_1\d_2\d_3\d_4\d_5\d_6\d_7\d_8\d_9 =
\s_4\s_7\a_2\b_7\s_5\s_8\a_8 \b_4\s_3\s_6\a_5\b_1$$
 on the nine-holed torus induces the word
$(\a^3\b_\a )^3=1$ in the mapping class group $\Gamma_1$ which gives
us an elliptic Lefschetz fibration on the elliptic surface
$E(1)=$\,\cp $\#9$  \cpb, where $\b_\a = \a\b\overline{\a}$ (which
is indeed  a right-handed Dehn twist). Note that we cyclicly
permuted the curves in the equation we derived in
Section~\ref{nine-star} to obtain  the  relation above.

\begin{figure}[ht]

  \begin{center}

     \includegraphics[width=7cm]{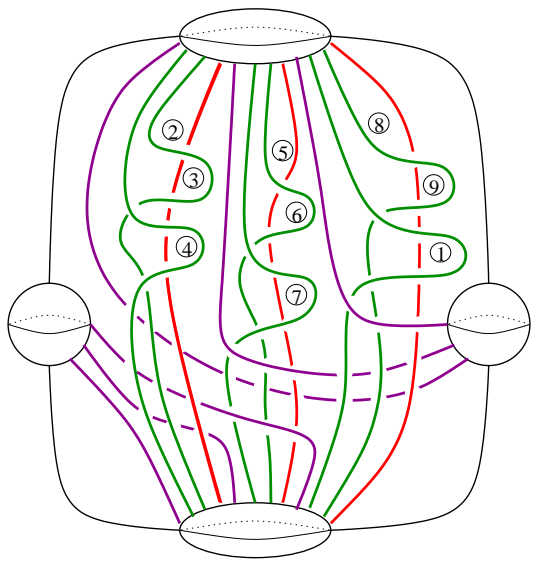}

   \caption{An elliptic Lefschetz fibration $E(1) \to S^2$ with nine disjoint sections.
   Note that we just indicated
   the intersection points of sections with the regular fiber by the encircled numbers
   corresponding to the boundary components $\d_1, \ldots, \d_9$ of the nine-holed torus.}

  \label{9-section}

    \end{center}

  \end{figure}

\vspace{1ex}

Finally, for $4\leq k \leq 9$, by taking the $n$-th power of our
relation for the $k$-holed torus we can find $k$ disjoint sections
of the corresponding elliptic fibration on the elliptic surface
$E(n)$ for any $n\geq 1$.

\section{Appendix}

In this appendix, we outline how to prove the lantern relation in Section~\ref{lant}
and the one-holed torus relation in Section~\ref{one-holedtorus}. These are well known,
but we include them here to have a complete presentation.

It is a classical fact that a self-diffeomorphism of a $2$-disc which is identity on the boundary
is isotopic to the identity through such diffeomorphisms. It follows that in order to show that
two self-diffeomorphisms of a surface are isotopic to each other, it is sufficient to prove that
their action on a set of curves whose complement is a disc are the same up to isotopy.

\begin{figure}[ht]

  \begin{center}

     \includegraphics{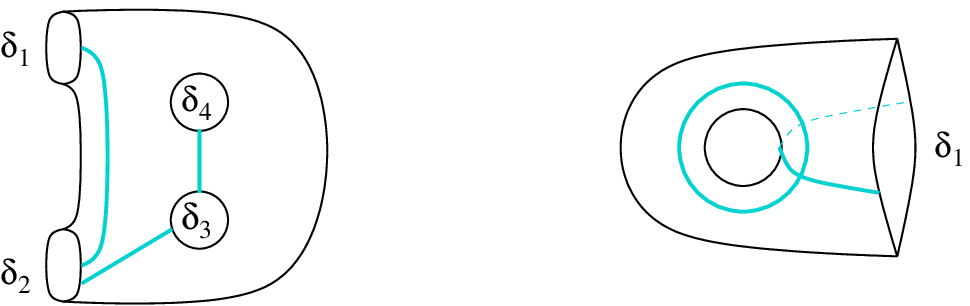}

   \caption{ }

  \label{rel}

    \end{center}

  \end{figure}

After this remark, the lantern relation $\d_1\d_2\d_3 \d_4 = \g \s \a$  is proved by looking at the action
of the diffeomorphisms on the both sides of the equality on the curves shown on the left in Figure~\ref{rel}.
Similarly, the one-holed torus relation $(\a \b)^6=\d_1$ is proved by looking at the action
of the diffeomorphisms on the both sides of the equality on the curves shown on the right in Figure~\ref{rel}.

\end{document}